\documentclass[11pt]{amsart}
\usepackage{graphicx}
\usepackage{amssymb}
\usepackage{amsmath}
\usepackage{amsthm}
\usepackage{amscd}
\usepackage[all,2cell,dvips]{xy}
\UseAllTwocells \SilentMatrices \oddsidemargin=10pt
\newtheorem{thm}{Theorem}[section]

\newtheorem{cor}[thm]{Corollary}
\newtheorem{lem}[thm]{Lemma}
\newtheorem{exm}[thm]{Example}

\theoremstyle{definition}

\theoremstyle{remark}
\newtheorem{rem}[thm]{\bf Remark}
\numberwithin{equation}{section}

\makeatletter

\newcommand{\Rmnum}[1]{\expandafter\@slowromancap\romannumeral #1@}
\makeatother

\begin{document}
\title[Singularity categories, Schur Functors and Triangular Matrix Rings]{Singularity categories,
Schur Functors and Triangular Matrix Rings}
\author[Xiao-Wu Chen]{Xiao-Wu Chen}
\thanks{Supported in part by the National Nature Science Foundation of China (Grant
No.s 10501041, 10601052)}
\thanks{E-mail: xwchen$\symbol{64}$mail.ustc.edu.cn}
\maketitle
\date{}%
\dedicatory{}%
\commby{}%

\begin{center}
Department of Mathematics\\
 University of Science and Technology of
China\\
Hefei 230026, P. R. China
\end{center}

\vskip 5pt

 \centerline{\footnotesize  Dedicated to Professor Freddy Van Oystaeyen
on the occasion of his sixtieth birthday}

\begin{abstract}
We study certain Schur functors which preserve singularity
categories of rings and we apply them to study the singularity
category of triangular matrix rings. In particular, combining these
results with Buchweitz-Happel's theorem, we can describe singularity
categories of certain non-Gorenstein rings via the stable category
of maximal Cohen-Macaulay modules. Three concrete examples of
finite-dimensional algebras with the same singularity category are
discussed.
\end{abstract}

\section{Introduction}

Singularity category is an important invariant for rings of infinite
global dimension and for singular varieties (\cite{Buc,Ha2}).
Recently, Orlov rediscovers the notion of singularity categories
(\cite{O1,O2,O3}) in his study of B-branes on Landau-Ginzburg models
in the framework of Homological Mirror Symmetry Conjecture (compare
\cite{KL}). Orlov shows that the category of B-branes on
Landau-Ginzburg models (proposed by Kontsevich) is equivalent to the
products of some singularity categories (\cite{O1}, Theorem 3.9 and
Corollary 3.10); he shows that the singularity category of algebraic
varieties enjoys the local property (\cite{O1}, Proposition 1.14).
Meanwhile, the singularity category of non-commutative rings and
algebras is also a very active topic, and there are extensive
references on it, see the introduction and references of \cite{CZ}.
It is known due to Buchweitz \cite{Buc} and independently Happel
\cite{Ha2} that the singularity categories of a Gorenstein ring can
be characterized by the stable category of its maximal
Cohen-Macaulay modules. There are a number of important consequences
of this result, for example, from it we know that if two Gorenstein
rings are derived equivalent, then their stable categories of
maximal Cohen-Macaulay modules are triangle-equivalent; that the
singularity category of a Gorenstein artin algebra is Krull-Schmidt
and has Auslander-Reiten triangles (compare \cite{Ha2} and
\cite{AR}). However, for non-Gorenstein rings and algebras, very
little is known about their singularity categories.

The aim of this paper is twofold: (i) We prove a local property for
the singularity categories of non-commutative rings, using Schur
functors. See Theorem 2.1 and compare Orlov's result (\cite{O1},
Proposition 1.14 or \cite{O3}, Proposition 1.3); (ii) We apply
Theorem 2.1 and Buchweitz-Happel's theorem to characterize the
singularity categories of certain (upper) triangular matrix
(non-Gorenstein) rings via the stable category of some maximal
Cohen-Macaulay modules. See Theorem 4.1(1) and Corollary 4.2. We
give an easy criterion (Theorem 3.3) on verifying when an upper
triangular matrix ring is Gorenstein, from which one sees that the
rings in Corollary 4.2 may be non-Gorenstein. The singularity
categories of some concrete examples of finite-dimensional
(non-Gorenstein) algebras are given explicitly in the last section.

\section{Schur functors preserving singularity categories}

\subsection{}

Throughout, $R$ will be a left-noetherian ring with a unit. Denote
by $R\mbox{-mod}$ the category of finitely-generated left
$R$-modules and $R\mbox{-proj}$ its full subcategory of
finitely-generated projective modules. Recall that
$D^b(R\mbox{-mod})$ is the bounded derived category of
$R\mbox{-mod}$, and $K^b(R\mbox{-proj})$ is the bounded homotopy
category of $R\mbox{-proj}$. View $K^b(R\mbox{-proj})$ as a thick
triangulated subcategory of $D^b(R\mbox{-mod})$. The singularity
category (\cite{O1,O2,O3}) of $R$ is defined to be the following
Verdier quotient category
\begin{align*}
D_{\rm sg}(R):= \; D^b(R\mbox{-mod})/ {K^b(R\mbox{-proj})}.
\end{align*}
The singularity category reflects certain singularity of the ring
$R$.

It is known that if two left-noetherian rings are derived
equivalent, then they have the same singularity category. However
the converse is not true. The main theorem in this section is to
provide certain equivalence of singularity categories via Schur
functors.

Let $e$ be an idempotent of $R$. The Schur functor (\cite{Gr},
Chapter 6) is defined to be
\begin{align*}
S_e=eR \otimes_R - \; : \; R\mbox{-mod}\longrightarrow
eRe\mbox{-mod}
\end{align*}
where $eR$ is viewed a natural $eRe$-$R$-bimodule via the
multiplication map. We denote the kernel of $S_e$ by
$\mathcal{B}_e$. Then it is not hard to see that $\mathcal{B}_e$ is
a full abelian subcategory, and an $R$-module $M\in \mathcal{B}_e$
if and only if $eM=0$, and if and only if $(1-e)M=M$.

To state the main theorem, we need to introduce some notions which
are somehow inspired by the notion of regular and singular points in
algebraic geometry. An idempotent $e\in R$ is said to be
\emph{regular}, if for any module $M\in \mathcal{B}_{1-e}$, ${\rm
proj.dim}\; _RM< \infty$, where we denote by ${\rm proj.dim}\; _RM$
the projective dimension of $M$. If $e$ is not regular, we say that
$e$ is \emph{singular}. The idempotent $e$ is said to be
\emph{singularly-complete}, if $1-e$ is regular. Note that the
properties defined above are invariant under conjugations. Let us
remark that one may compare them with the notions in \cite{DEN},
2.3.

\vskip 10pt

 Our main result is

\begin{thm}
Let $R$ be a left-noetherian ring, $e$ its idempotent. Assume that
$e$ is singularly-complete and ${\rm proj.dim}\; _{eRe}eR < \infty$.
Then the Schur functor $S_e$ induces an equivalence of triangulated
categories $D_{\rm sg}(R) \simeq D_{\rm sg}(eRe)$.
\end{thm}

Let us remark that the theorem is inspired by a result of Orlov on
the local property of singularity categories of algebraic varieties.
Let $\mathbb{X}$ be an algebraic variety. Denote by $D_{\rm
sg}(\mathbb{X})$ the singularity category of $\mathbb{X}$ which is
defined as the Verdier quotient category of the bounded derived
category $D^b( {\rm coh}(\mathbb{X}))$ of coherent sheaves with
respect to its full triangulated category of perfect complexes ${\rm
perf}(\mathbb{X})$. In \cite{O1}, Proposition 1.14, Orlov shows the
following local property of singularity categories: that if
$\mathbb{X}'\subseteq \mathbb{X}$ is an open subvariety containing
the singular locus, then we have a natural equivalence of
triangulated categories $D_{\rm sg}(\mathbb{X})\simeq D_{\rm
sg}(\mathbb{X}')$. Thus in our situation, $eRe$ is viewed as an open
``subvariety'' of $R$, and the idempotent $e$ is singularly-complete
is somehow similar to saying that $eRe$ contains all the
``singularity'' of $R$. Therefore, our result can be regarded as a
(possibly naive) version of local property of non-commutative
singularity categories.

\subsection{}
Let us begin with an easy lemma. Let $R$ be a left-noetherian ring,
$e$ its idempotent. Then it is not hard to see that $eRe$ is also
left-noetherian. Recall the category $\mathcal{B}_e =\{M\in
R\mbox{-mod}\;| \; eM=0\}$. Let $\mathcal{N}_e$ be the full
subcategory of $D^b(R\mbox{-mod})$ consisting of complex $X^\bullet$
with its cohomology groups $H^n(X^\bullet)$ lying in
$\mathcal{B}_e$. It is a triangulated subcategory.

\begin{lem}
Use above notation. Then the Schur functor $S_e$ induces a natural
equivalence of triangulated categories $D^b(R\mbox{\rm
-mod})/{\mathcal{N}_e}\simeq D^b(eRe\mbox{\rm -mod})$.
\end{lem}

\noindent {\bf Proof.}\quad Note that the Schur functor $S_e$ is
exact and recall that the subcategory $\mathcal{B}_e$ is the kernel
of $S_e$, one sees that $\mathcal{B}_e$ is a Serre subcategory, and
it is well-known that the functor $S_e$ induces an equivalence of
abelian categories
$$R\mbox{-mod}/\mathcal{B}_e \simeq
eRe\mbox{-mod}.$$
Now the result follows immediately from a
fundamental result by Miyachi (\cite{Mi}, Theorem 3.2; also see the
appendix in \cite{BO}): for any abelian category $\mathcal{A}$ and
its Serre subcategory $\mathcal{B}$, we have a natural equivalence
of triangulated categories
$D^b(\mathcal{A})/{D^b(\mathcal{A})_\mathcal{B}}\simeq
D^b(\mathcal{A}/\mathcal{B})$, where
$D^b(\mathcal{A})_\mathcal{B}:=\{X^\bullet \in D^b(\mathcal{A})\; |
\; H^n(X^\bullet)\in \mathcal{B},\; n\in \mathbb{Z}\}$. \hfill
$\blacksquare$

\vskip 10pt

Let us recall some notions. Let $\mathcal{C}$ be a triangulated
category, $[1]$ its shift functor. Let $S\subseteq \mathcal{C}$ be a
subset. The smallest triangulated subcategory of $\mathcal{C}$
containing $S$ is denoted by $\langle S \rangle$, and it is said to
be generated by $S$. In fact, objects in $\langle S\rangle$ are
obtained by iterated extensions of objects from $\bigcup_{n\in
\mathbb{Z}} S[n]$, see \cite{Ha1}, p.70. For example,
$K^b(R\mbox{-proj})$ is generated by $R\mbox{-proj}$. Note that the
category $\mathcal{N}_e$ defined above is generated by
$\mathcal{B}_e$ (for example, by \cite{Har}, Lemma 7.2(4)).

\vskip 10pt

 \noindent {\bf Proof of Theorem 2.1.}\quad Since $e$ is
 singularly-complete, then every module in $\mathcal{B}_{e}$ has
 finite projective dimension. Hence inside $D^b(R\mbox{-mod})$, we
 have $\mathcal{B}_e\subseteq K^b(R\mbox{-proj})$. Since
 $\mathcal{N}_e$ is generated by $\mathcal{B}_e$, we get $\mathcal{N}_e\subseteq
 K^b(R\mbox{-proj})$. Consequently, we have
 \begin{align*}
 D_{\rm sg}(R)= D^b(R\mbox{-mod})/{K^b(R\mbox{-proj})} \simeq
 (D^b(R\mbox{-mod})/{\mathcal{N}_e})/(K^b(R\mbox{-proj})/\mathcal{N}_e).
 \end{align*}
By Lemma 2.2, the Schur functor $S_e$ induces a natural equivalence
$$\bar{S_e}: D^b(R\mbox{-mod})/{\mathcal{N}_e}\simeq
D^b(eRe\mbox{-mod}).$$ Therefore it suffices to show that the
essential image of $K^b(R\mbox{-proj})/{\mathcal{N}_e}$ under
$\bar{S_e}$ is exactly $K^b(eRe\mbox{-proj})$.

 To see this, denote the essential image by $\mathcal{M}$. Since
 $K^b(R\mbox{-proj})$ is generated by $R\mbox{-proj}$, hence
 $\mathcal{M}$ is generated by $S_e(R\mbox{-proj})$. By the
 assumption, $S_e(R)=eR$ has finite projective dimension over $eRe$,
 hence for every projective $R$-module $P$, $S_e(P)$ has finite
 projective dimension, in other words, we have $S_e(R\mbox{-proj})\subseteq
 K^b(eRe\mbox{-proj})$, and therefore $\mathcal{M}\subseteq
 K^b(eRe\mbox{-proj})$. On the other hand, note that $S_e$ induces
 an equivalence of categories
 $${\rm add}\; Re \simeq eRe\mbox{-proj},$$
 where ${Re}$
 is the projective left $R$-module determined by $e$ and ${\rm add}\; Re$ is the full subcategory
 of $R\mbox{-mod}$ consisting of all the direct summands of sums of finite copies of $Re$. Hence we know that
 $S_e(R\mbox{-proj})$ contains a set of generators for
 $K^b(eRe\mbox{-proj})$, and thus we obtain that $\mathcal{M}$ contains
 $K^b(eRe\mbox{-proj})$. Thus we are done. \hfill $\blacksquare$

\section{Triangular Matrix Gorenstein Rings}

In this section, we will study triangular matrix rings. The main
result states an easy criterion on when an upper triangular matrix
ring is Gorenstein.

\subsection{}  Recall some facts on (upper) triangular matrix rings
(compare \cite{ARS}, III.2). Let $R$ and $S$ be any rings,
$M={_RM_S}$ an $R$-$S$-bimodule. We study the corresponding upper
triangular matrix ring $T= \begin{pmatrix}R & M
\\
0 & S\end{pmatrix}$.

Recall the description of left $T$-modules via column vectors. Given
a left $R$-module $_RX$ and a left $S$-module $_SY$, and an
$R$-module morphism $\phi: M \otimes_S Y \longrightarrow X$, we
define the left $T$-module structure on $\begin{pmatrix} X \\
Y\end{pmatrix}$ by the following identity
\begin{align*}
\begin{pmatrix} r & m \\
                0 & s\end{pmatrix}\;
 \begin{pmatrix} x\\ y \end{pmatrix} :=
\begin{pmatrix} rx+\phi(m\otimes y) \\ sy \end{pmatrix}.
\end{align*}
It is not hard to check that every $T$-module arises in this way
(compare \cite{ARS}, III.2, Proposition 2.1).

The following lemma is well-known, and it could be checked directly
(compare \cite{ARS}, III, Proposition 2.3 and 2.5(c)).

\begin{lem}
Use above notation.  \\
(1).\quad We have  ${\rm proj.dim}\; \begin{pmatrix} X \\ 0
\end{pmatrix} = {\rm proj.dim} \; {_RX}, \; {\rm inj.dim}\;
\begin{pmatrix} 0 \\ Y
\end{pmatrix}={\rm inj.dim}\; {_SY}$.\\
(2).\quad  For any $R$-module $_RX'$, we have a natural isomorphism
\begin{align*}
{\rm Hom}_T(\begin{pmatrix} X\\ Y \end{pmatrix},\; \begin{pmatrix}
X' \\ {\rm Hom}_R(M, X')\end{pmatrix})\simeq {\rm Hom}_R(X, X'),
\end{align*}
where $\begin{pmatrix} X' \\ {\rm Hom}_R(M, X')\end{pmatrix}$
becomes a left $T$-module via the natural evaluation map $M\otimes_S
{\rm Hom}_R(M, X')\longrightarrow X'$. In particular,
$\begin{pmatrix} X'
\\ {\rm Hom}_R(M, X')\end{pmatrix}$ is an injective $T$-module if
and only if $_RX'$ is injective.
\end{lem}

Dually, we have the description of right $T$-modules via row
vectors. Precisely, given a right $R$-module $X_R$ and a right
$S$-module $Y_S$, and a right $S$-module morphism $\psi: X\otimes_R
M \longrightarrow Y$, then the space $\begin{pmatrix} X & Y
\end{pmatrix}$ carries a right $T$-module structure via the
following identity
\begin{align*}
\begin{pmatrix}x & y \end{pmatrix} \; \begin{pmatrix} r & m \\
                                                      0 &  s
                                                      \end{pmatrix}:=
                                                      \begin{pmatrix}
                                                      xr & \psi(x\otimes
                                                      m)+ys
                                                      \end{pmatrix}.
\end{align*}

Dual to Lemma 3.1, we have

\begin{lem}
Use above notation.  \\
(1).\quad We have  ${\rm inj.dim}\; \begin{pmatrix} X & 0
\end{pmatrix} = {\rm inj.dim} \; {X_R}, \; {\rm proj.dim}\;
\begin{pmatrix} 0 & Y
\end{pmatrix} = {\rm proj.dim}\; {Y_S}$.\\
(2).\quad  For any right $S$-module $Y'_S$, we have a natural
isomorphism
\begin{align*}
{\rm Hom}_T(\begin{pmatrix} X& Y \end{pmatrix},\; \begin{pmatrix}
{\rm Hom}_S(M, Y') &Y'\end{pmatrix})\simeq {\rm Hom}_R(Y, Y'),
\end{align*}
where $\begin{pmatrix} {\rm Hom}_S(M, Y') & Y'\end{pmatrix}$ becomes
a right $T$-module via the natural evaluation map ${\rm Hom}_S(M,
Y')\otimes_R M\longrightarrow Y'$. In particular, $\begin{pmatrix}
{\rm Hom}_S(M, Y')& Y'\end{pmatrix}$ is an injective $T$-module if
and only if $Y'_S$ is injective.
\end{lem}
\vskip 5pt

Similarly, we may consider lower triangular matrix rings. Let $R$,
$S$ be rings, $M={_RM_S}$ be a bimodule. Then we have the lower
triangular matrix ring $T'= \begin{pmatrix}S & 0
\\
M & R\end{pmatrix}$. Note that the opposite ring $T'^{\rm op}$ is an
upper triangular matrix ring, in fact, $T'^{\rm op} =
\begin{pmatrix}S^{\rm op} & M
\\
0 & R^{\rm op}\end{pmatrix}$, where $M$ is viewed as an $S^{\rm
op}$-$R^{\rm op}$-bimodule. Hence, one can deduce easily the
corresponding results for lower triangular matrix rings from Lemma
3.1 and 3.2. We will quote these results without writing them down
explicitly.

\subsection{}

Recall that a ring $R$ is said to be Gorenstein, if $R$ is two-sided
noetherian and the regular module $R$ has finite injective dimension
both as left and right modules \cite{EJ}. An artin algebra which is
Gorenstein is called a Gorenstein artin algebra (\cite{Ha2} or
\cite{AR}). It is shown by Zaks (\cite{Z}, Lemma A) that for
Gorenstein ring $R$, we have ${\rm inj.dim}\; _RR={\rm inj.dim}\;
R_R$, while this integer will be denoted by ${\rm G.dim}\;R$.

\vskip 5pt

We have the main result in this section.

\begin{thm} With above notion. Assume that both $R$ and $S$ are
Gorenstein rings, $M={_RM_S}$ an $R$-$S$-bimodule. Then the upper
triangular matrix ring
$\begin{pmatrix} R & M \\
                                               0 & S
\end{pmatrix}$ is Gorenstein if and only if both $_RM$ and $M_S$ are
finitely-generated, ${\rm proj.dim}\; _RM< \infty$ and ${\rm
proj.dim}\; M_S < \infty$.
\end{thm}

Before giving the proof, we note the following basic fact.

\begin{lem}{\rm (\cite{EJ}, Proposition 9.1.7)}
Let $R$ be a Gorenstein ring, $M={_RM}$ a left $R$-module. Then $M$
has finite projective dimension if and only if $M$ has finite
injective dimension.
\end{lem}

\vskip 10pt

\noindent {\bf Proof of Theorem 3.3.}\quad  Denote by $T$ the upper
triangular matrix ring in our consideration. The ``only if'' part is
easy. Assume that $T$ is Gorenstein. Consider the following exact
sequence of left $T$-modules:

\begin{align}
0\longrightarrow \begin{pmatrix} M \\ 0 \end{pmatrix}
\longrightarrow \begin{pmatrix} M \\ S \end{pmatrix} \longrightarrow
\begin{pmatrix} 0 \\ S\end{pmatrix} \longrightarrow 0.
\end{align}
Note that the middle term is a principal module (i.e., a cyclic
projective module), in particular, it is noetherian. Hence the
$T$-module $\begin{pmatrix} M \\ 0 \end{pmatrix}$ is noetherian, and
it follows immediately that $_RM$ is noetherian. Moreover, since
$_SS$ has finite injective dimension, and by Lemma 3.1(1), the last
term has finite injective dimension, and then by Lemma 3.4, it has
finite projective dimension. Now it follows that ${\rm
proj.dim}\;\begin{pmatrix} M
\\ 0 \end{pmatrix}$ is finite. By Lemma 3.1(1) again, we get ${\rm proj.dim}\;
{_RM}<\infty$. Similarly, one can show that $M_S$ is noetherian and
${\rm proj.dim} M_S<\infty$.

\vskip 5pt

 Next we show the ``if'' part. Assume that both $_RM$ and
$M_S$ are finitely-generated and hence noetherian,  and have finite
projective dimension. First consider the following exact sequence of
left (or right) $T$-modules
\begin{align}
0 \longrightarrow \begin{pmatrix} 0 & M \\ 0 & 0 \end{pmatrix}
\longrightarrow T \longrightarrow \begin{pmatrix} R & 0 \\ 0 & S
\end{pmatrix} \longrightarrow 0.
\end{align}
From the assumption, we know the first term is a left noetherian
$T$-module; because of the noetherianness of $R$ and $S$, the last
term viewed as a left $T$-module is also noetherian. Hence $_TT$ is
noetherian, that is, $T$ is left-noetherian. Similarly, $T$ is
right-noetherian, and thus $T$ is a two-sided noetherian ring. What
is left to show is that ${\rm inj.dim} \; {_TT}< \infty$ and ${\rm
inj.dim}\; {T_T}<\infty$. We only show that ${\rm inj.dim}\;
{_TT}<\infty$, and the other can be proven similarly (using Lemma
3.2).

To show that ${\rm inj.dim}\; {_TT}< \infty$, first note we have a
decomposition of left $T$-mdoules $_TT= \begin{pmatrix}R \\
0\end{pmatrix} \oplus \begin{pmatrix} M \\ S \end{pmatrix}$. We
claim that it suffices to show that ${\rm inj.dim}\; \begin{pmatrix}R \\
0\end{pmatrix}< \infty$. In fact, since $_RM$ has finite projective
dimension, we have an exact sequence of left $R$-modules
$$0 \longrightarrow P^m
\longrightarrow P^{m-1} \longrightarrow \cdots \longrightarrow
P^0\longrightarrow M\longrightarrow 0,$$
 where each $P^j$ is a
finitely-generated projective $R$-module. Since ${\rm inj.dim}
\begin{pmatrix} R \\ 0\end{pmatrix}< \infty$, we know that each
$T$-module $\begin{pmatrix}P^j \\
0\end{pmatrix}$ has finite injective dimension, and note the
following natural exact sequence of $T$-modules
$$0 \longrightarrow \begin{pmatrix}P^m \\ 0\end{pmatrix}
\longrightarrow \begin{pmatrix}P^{m-1} \\ 0\end{pmatrix}
\longrightarrow \cdots \longrightarrow
\begin{pmatrix}P^0\\ 0 \end{pmatrix}\longrightarrow \begin{pmatrix}M\\ 0 \end{pmatrix} \longrightarrow 0,$$
thus we obtain that ${\rm inj.dim}\; \begin{pmatrix}M \\
0\end{pmatrix}< \infty.$ Now by (3.1) and note that ${\rm inj.dim}\;
\begin{pmatrix}0 \\ S\end{pmatrix}< \infty$, we get that ${\rm inj.dim}\; \begin{pmatrix}M \\ S\end{pmatrix}<
\infty$. Thus we are done with ${\rm inj.dim}\; {_TT}<\infty$.

To prove ${\rm inj.dim}\; \begin{pmatrix}R \\
0\end{pmatrix}< \infty$, since $_RR$ has finite injective dimension,
we may take its finite injective resolution. Applying the same
argument as above, one deduces that it suffices to show that ${\rm
inj.dim}\;
\begin{pmatrix} I\\0 \end{pmatrix}<\infty$ for each injective
$R$-module $I$. Consider the following natural exact sequence of
$T$-modules
\begin{align*}
0 \longrightarrow \begin{pmatrix}I \\ 0\end{pmatrix} \longrightarrow
\begin{pmatrix}I \\ {\rm Hom}_R(M, I)\end{pmatrix} \longrightarrow
\begin{pmatrix}0 \\ {\rm Hom}_R(M, I)\end{pmatrix} \longrightarrow
0.
\end{align*}
By Lemma 3.1(2), the middle term is injective. Using ${\rm
proj.dim}\; _RM<\infty$, it is easy to show by the Hom-tensor
adjoint that ${\rm inj.dim}\; _S {\rm Hom}_R(M, I)< \infty$. By
Lemma 3.1(1), the last term has finite injective dimension, and thus
so does the first term. This completes the proof. \hfill
$\blacksquare$

\begin{rem}
From the proof above and using dimension-shift if necessary, it is
not hard to see that: in the situation of Theorem 3.3, we have
\begin{align*}
{\rm max}\{{\rm G.dim}\; R, {\rm G.dim}\; S\} \leq {\rm G.dim}\;
\begin{pmatrix} R & M \\ 0 & S\end{pmatrix} \leq {\rm G.dim}\; R +
{\rm G.dim}\; S +1.
\end{align*}
\end{rem}

\section{Applications and Examples}

This section is devoted to applying the above results to the
singularity categories of certain (non-Gorenstein) rings and
algebras. Three concrete examples are included

\subsection{}

Recall that a ring $R$ is said to be regular, if $R$ is two-sided
noetherian and $R$ has finite global dimension on both sides. The
following application of Theorem 2.1 is our main result.

\begin{thm}
Let $R$ be a left-noetherian ring with finite left global dimension,
$S$ a left-noetherian ring. \\
(1)\quad Let $M={_RM_S}$  be a bimodule such that $_RM$ is
finitely-generated. Then we have a natural equivalence of
triangulated categories $$D_{\rm sg}(\begin{pmatrix} R & M \\
                                                   0 & S \end{pmatrix})\simeq D_{\rm
                                                   sg}(S). $$
(2)\quad Let $N={_SN_R}$ be a bimodule such that both $_SN$ and
$N_R$ are finitely-generated and  ${\rm proj.dim}\; _SN< \infty$.
Assume furthur that $R$ is regular and $S$ is Gorenstein. Then we
have a natural
equivalence of triangulated categories $$D_{\rm sg}(\begin{pmatrix} R & 0 \\
                                                   N & S \end{pmatrix})\simeq D_{\rm
                                                   sg}(S). $$ \end{thm}

\noindent {\bf Proof.}\quad (1)\quad Denote by $T$ the upper
triangular matrix ring in our consideration. By (3.2), one deduces
easily that $T$ is left-noetherian and thus its singularity category
is well-defined.
Set $e:= \begin{pmatrix} 0 & 0 \\
0 & 1\end{pmatrix}$, and thus $eTe\simeq S$. In this case, we have
$\mathcal{B}_{e}= \{
\begin{pmatrix} X \\ 0 \end{pmatrix}\; | \; X={_RX} \mbox{ any } R
\mbox{-module}\}$. Since $R$ has finite left global dimension, then
${\rm proj.dim}\; {_RX}< \infty
$, and thus by Lemma 3.1(1), we get ${\rm proj.dim}\; \begin{pmatrix} X \\
0
\end{pmatrix}< \infty$, therefore $e$ is singularly-complete. It is
easy to see, as an $eTe$-module, $eT=eTe$, and hence ${\rm
proj.dim}\; {_{eTe}eT}=0$. Thus the conditions of Theorem 2.1 are
fulfilled, and the result follows.

\vskip 5pt
 (2) \quad Denote by $T$ the lower triangular matrix ring in
this consideration. As above, it is not hard to see that the ring
$T$ is left-noetherian and thus $D_{\rm sg}(T)$ is defined. Note
that $T^{\rm op}$ is an upper triangular matrix ring (see 3.1), and
by Theorem 3.3, $T^{\rm op}$ and thus $T$ is Gorenstein. We still
denote left $T$-modules by column vectors. As above,
set $e:= \begin{pmatrix} 0 & 0 \\
0 & 1\end{pmatrix}$, and $\mathcal{B}_{1-e}= \{
\begin{pmatrix} X \\ 0 \end{pmatrix}\; | \; X={_RX} \mbox{ any } R
\mbox{-module}\}$. View left $T$-module as right $T^{\rm
op}$-modules. Using Lemma 3.2 (1), we deduce that ${\rm inj.dim}\;
\begin{pmatrix} X \\ 0 \end{pmatrix}={\rm inj.dim}\; _RX< \infty$.
Since $T$ is Gorenstein, and then by Lemma 3.4, we get ${\rm
proj.dim}\; \begin{pmatrix} X \\ 0 \end{pmatrix}< \infty$. Hence the
idempotent $e$ is singularly-complete.

 Next we show that ${\rm proj.dim}\; _{eTe}{eT}<\infty$. Since $eT=eTe\oplus
 eT(1-e)$, and hence it suffices to show that ${\rm proj.dim}\; _{eTe}{eT(1-e)}<
 \infty$. Note that $eTe=S$ and, viewed as a left $eTe$-module
 $eT(1-e)\simeq N$, by the assumption, ${\rm proj.dim} \; {_SN}< \infty$, and thus
 we obtain that ${\rm proj.dim}\; _{eTe}{eT}<\infty$. Therefore
 the conditions of Theorem 2.1 are fulfilled, and the result
 follows. \hfill $\blacksquare$

\vskip 10pt

Theorem 4.1(1) allows us to describe the singularity categories of
some non-Gorenstein rings as the stable category of certain maximal
Cohen-Macaulay modules. For this end, let us recall a result by
Buchweitz (\cite{Buc}, Theorem 4.4.1) and independently by Happel
(\cite{Ha2}, Theorem 4.6; compare \cite{CZ}, Theorem 2.5). Let $R$
be a Gorenstein ring. Denote by
$${\rm MCM(R)}:=\{M \in R\mbox{-mod}\; |
\; {\rm Ext}^i_R(M, R)=0, \; i\geq 1\} $$
 the category of maximal
Cohen-Macaulay modules. It is a Frobenius category with (relative)
projective-injective objects exactly contained in $R\mbox{-proj}$
(compare \cite{CZ}, 2.1). Denote by $\underline{\rm MCM}(R)$ its
stable category, which is a triangulated category (\cite{Ha1},
p.16). Then Buchweitz-Happel's theorem says that there is  an
equivalence of triangulated categories $D_{\rm sg}(R)\simeq
\underline{\rm MCM}(R)$. This generalizes a result of Rickard
\cite{Ric}, which says that for a self-injective algebra, its
singularity category is triangle-equivalent to the stable category
of its module category. However for non-Gorenstein rings and
algebras, we know little about their singularity categories.

The following result is a direct consequence of Theorem 4.1(1) and
Buchweitz-Happel's theorem.

\begin{cor}
 Let $R$ be a regular ring, $S$ a Gorenstein ring, $M= {_RM_S}$ a bimodule which is finitely-generated on $R$.
 Then we have an equivalence of triangulated categories
 $$D_{\rm sg}(\begin{pmatrix} R & M \\
                                               0 & S
\end{pmatrix})\simeq \underline{\rm MCM}(S).$$
\end{cor}
 Note that by Theorem 3.3 the above upper triangular matrix ring will be non-Gorenstein,
 provided that  ${\rm proj.dim}\; M_S=\infty$.

\subsection{}
 In this subsection, we will study three examples of
 finite-dimensional
 algebras, which share the same singularity category. Let us remark that these
 examples can be easily generalized. $K$  will be a field.

\begin{exm}
(1).\quad Let $A$ be the  $K$-algebra given by the following quiver
and relations
\[
\xymatrix@C=50pt @R=30pt{\cdot_1 \ar@(l,u)^\alpha \ar@/^/[r]^\beta &
\cdot_2 \ar@/^/[l]^\gamma }\quad \quad
\alpha^2=\gamma\beta=0=\beta\alpha.\]

 Here we write the
concatenation of paths from the left to the right. Set $e=e_1$. Note
that the second simple module $S_2$ has finite projective dimension,
and every module in $\mathcal{B}_{1-e}$ is obtained by iterated
extensions of $S_2$, and thus of finite projective dimension. Hence
the idempotent $e$ is a singularly-complete idempotent. It is not
hard to see that $eAe\simeq K[x]/{(x^2)}$, and $eA$ is a free left
$eAe$-module. Hence by Theorem 2.1, we have an equivalence of
triangulated categories
 $$D_{\rm sg}(A)\simeq D_{\rm sg}(K[x]/{(x^2)}).$$

 Since $K[x]/{(x^2)}$ is Frobenius, then by Rickard's theorem (\cite{Ric}, Theorem 2.1),
 we have a triangle-equivalence: $D_{\rm
sg}(K[x]/{(x^2)})\simeq K[x]/{(x^2)}\mbox{-}\underline{\rm mod}$.

Recall that every semisimple abelian category (for example, the
category $K\mbox{\rm -mod}$ of finite-dimensional $K$-spaces), in a
unique way, becomes a triangulated category with the identity
functor being the shift functor. Then it is not hard to see that
there is a triangle-equivalence $K[x]/{(x^2)}\mbox{-}\underline{\rm
mod}\simeq K\mbox{\rm -mod}$. Hence we finally get a
triangle-equivalence
$$D_{\rm sg}(A) \simeq K\mbox{\rm -mod}.$$

Let us remark that the algebra $A$ is not Gorenstein since ${\rm
proj.dim}\; I(2)=\infty$, where $I(2)$ is the injective hull of
$S_2$.

\vskip 10pt

 (2).\quad  Let  $A'$ be the $K$-algebra given  the following quiver and relations
\[
\xymatrix@C=50pt @R=30pt { \cdot_1  \ar@(l,u)^\alpha \ar[r]^\beta &
\cdot_2 }  \quad  \quad \alpha^2=0=\beta \alpha.\]

Then one may view $A'$ as an upper triangular matrix algebra:
$A'=\begin{pmatrix} e_2A'e_2 & e_2A'e_1 \\ 0 &
e_1A'e_1\end{pmatrix}$.  Note that $e_2A'e_2\simeq K$,
$e_1A'e_1\simeq K[x]/{(x^2)}$, and $e_2A'e_1$ is not a projective
$e_1A'e_1$-module. Since $e_1A'e_1$ is Frobenius, thus one infers
that $e_2A'e_1$, as a right $e_1A'e_1$-module, is of infinite
projective dimension. By Theorem 3.3, we deduce that $A'$ is not
Gorenstein. However by Corollary 4.2, it is not hard to see that
$D_{\rm sg}(A')\simeq K[x]/(x^2){\rm-\underline{mod}}$, and thus we
get a triangle-equivalence
$$D_{\rm sg}(A')\simeq K\mbox{\rm
-mod}.$$

\vskip 10pt

(3)\quad  Let  $A''$ be the $K$-algebra given  the following quiver
and relations
\[
\xymatrix@C=50pt @R=30pt { \cdot_1  \ar@(l,u)^\alpha  & \cdot_2
\ar@/^/[l]^\beta \ar@/_/[l]_\gamma } \quad  \quad \alpha^2=0. \]

Then the algebra $A''$ can be viewed a lower triangular matrix
algebra $A''=\begin{pmatrix} e_2A''e_2 & 0\\  e_1A''e_2 & e_1A''e_1
\end{pmatrix}$. Note that $e_2A''e_2\simeq K$, $e_1A''e_1\simeq
K[x]/{(x^2)}$, and $e_1A''e_2$, viewed as a left $e_1A''e_1$-module,
is free of rank 2. By Theorem 4.1(2), we have an equivalence of
triangulated categories $D_{\rm sg}(A'')\simeq D_{\rm
sg}(K[x]/{(x^2)})$, and thus by the argument above, we also have a
traingle-equivalence
$$D_{\rm sg}(A'')\simeq K\mbox{\rm
-mod}.$$ Note that by Theorem 3.3 (or rather the corresponding
result for lower triangular matrix rings), the algebra $A''$ is
Gorenstein.
\end{exm}

\vskip 5pt

 \noindent {\bf Acknowledgement}\quad The author would like to thank
 Prof. Ragnar-Olaf Buchweitz very much for sending him the beautiful
 paper \cite{Buc}, and to thank Prof. Pu Zhang for his interest in this work.

\bibliography{}

\end{document}